\documentclass[12pt]{article}%
\usepackage[applemac]{inputenc}
\usepackage{amsmath,amssymb}
\usepackage{amsthm}
\usepackage{amsfonts}
\usepackage{amsmath,amssymb}
\usepackage[applemac]{inputenc}
\usepackage{color}
\usepackage{color, colortbl}
\usepackage{xcolor}
\usepackage{amsmath}
\usepackage{amssymb}
\usepackage{graphicx}
\usepackage{tikz}
\usepackage{geometry}
\usepackage{mathtools}
\allowdisplaybreaks
\usepackage{float}           
\usepackage{graphicx}        
\usepackage{caption}         
\usepackage{subcaption}      
\usetikzlibrary{arrows}
\usepackage{amsfonts}
\usepackage{amsmath,amssymb}
\usepackage[applemac]{inputenc}
\usepackage{color}
\usepackage{amsmath}
\usepackage{amssymb}
\usepackage{graphicx}
\usepackage{tikz}
\usepackage[numbers]{natbib}  
\usepackage{needspace}

\newtheorem{theorem}{Theorem}[section]
\newtheorem{lemma}[theorem]{Lemma}
\newcommand{\E}{E}
\newtheorem{proposition}[theorem]{Proposition}
\newtheorem{corollary}[theorem]{Corollary}
\newtheorem{definition}[theorem]{Definition}

\newtheorem{example}[theorem]{Example}

\newtheorem{remark}[theorem]{Remark}

\usetikzlibrary{arrows,shapes,automata,petri}

\newcommand{\dproof}{\noindent {Proof.} \quad}
\newcommand{\fproof}{\hfill $\square$ \bigskip}

\numberwithin{equation}{section}

\definecolor{LightCyan}{rgb}{0.88,1,1}

\def\1B{\text{1\!\!I}}

\def\<{\langle}
\def\>{\rangle}

\def\F{\mathcal{F}}

\def\H{\mathcal{H}}

\def\S{\mathcal{S}}

\def\C{\mathcal{C}}
\def\E{\mathbb{E}}
\def\R{\mathbb{R}}
\def\P{\mathbb{P}}
\def\Q{\mathbb{Q}}

\def\J{\mathcal{J}}

\definecolor{strongred}{RGB}{200, 0, 0}

\definecolor{darkpurple}{HTML}{800080}

\definecolor{deepblue}{cmyk}{1,1,0,0.2}
\begin{document}

\title{Markov property  for systems driven by the  Brownian sheet}
\author{Nacira Agram$^{1}$ Bernt \O ksendal$^{2}$ Frank Proske$^{2}$  Olena Tymoshenko$^{2,3}$}
\date{29 September 2026}
\maketitle

\footnotetext[1]{Department of Mathematics, KTH Royal Institute of Technology and Digital Futures, 100 44, Stockholm, Sweden. 
Email: nacira@kth.se. Work supported by the Swedish Research Council grant (2020-04697).}

\footnotetext[2]{
Department of Mathematics, University of Oslo, Norway. 
Emails: oksendal@math.uio.no, proske@math.uio.no, olenaty@math.uio.no}
\footnotetext[3]{%
Department of Mathematical Analysis and Probability Theory, Igor Sikorsky Kyiv Polytechnic Institute, Ukraine. 
}

\begin{abstract}
We study stochastic systems driven by a Brownian sheet and formulate a
Markov property adapted to their two-parameter structure.  For time-space
homogeneous It\^o sheets, the future evolution from a point $(t,x)$ is not
determined by the single value $Y(t,x)$ alone, but also by the boundary data
along the horizontal and vertical edges issuing from $(t,x)$.  We therefore
show that the enlarged state consisting of the triple
$(Y,D_1Y,D_2Y)$ satisfies a natural Markov property with respect to the
past sigma-algebra generated by the two coordinate histories of the Brownian
sheet.  The result is illustrated by examples showing the strong martingale
property for the Brownian sheet and for simple Brownian sheet systems. We also recall a Brownian sheet Girsanov theorem and use it to treat drifted sheets.

\end{abstract}




\section{Introduction}

The Brownian sheet is the canonical two-parameter analogue of Brownian
motion.  It appears naturally in the study of random fields, stochastic
partial differential equations, and stochastic systems whose evolution is
indexed by time and space.  In contrast with one-parameter diffusion theory,
however, the two-parameter setting has a more delicate causal structure.  The
increment over a rectangle depends on information coming from two boundary
directions, and the usual one-dimensional formulation of the Markov property
does not directly apply.

The stochastic integration theory for Brownian sheets and martingale measures
was developed by Cairoli and Walsh \cite{CW} and Walsh \cite{Walsh}.
Markovian aspects of Brownian sheets and related multiparameter random fields
have been studied in several forms. In particular, germ-field Markov
properties for multiparameter processes were investigated by Mandrekar
\cite{Mandrekar}. Markov properties for two-parameter Gaussian processes were
studied by Nualart and Sanz \cite{NualartSanz}, Markov properties for certain
random fields were studied by Carnal and Walsh \cite{CarnalWal}, and sharp
Markov properties of the Brownian sheet and related processes were developed by
Dalang and Walsh \cite{DalangWal}. These works show that Brownian sheets and, more generally,
multiparameter random fields admit several Markov-type structures, and that
the appropriate notion of Markovianity depends strongly on the geometry of the
conditioning set.

Previous results establish Markov properties of the Brownian sheet and related multiparameter random fields relative to geometrically defined conditioning sets. Here we take a different viewpoint. Rather than studying the Brownian sheet itself as a random field, we consider nonlinear stochastic systems driven by the sheet and ask what information is required to restart their evolution from a point $(t,x)$. We show that the point value $Y(t,x)$ alone is not sufficient; one must also retain the boundary information along the two coordinate directions. This leads to an enlarged state consisting of the current value together with the two boundary profiles, for which we establish a natural Markov property with respect to the two-sided past sigma-algebra.

More precisely, we consider time-space homogeneous  It\^o sheets of the form
\begin{align*}
Y(t,x)&=y+
\int_{t_0}^t p(s)ds+
\int_{x_0}^x q(a)da
+\int_{t_0}^t\int_{x_0}^x b(Y(s,a))dsda  \\
&\quad+
\int_{t_0}^t\int_{x_0}^x \sigma(Y(s,a))B(dsda),
\qquad t\geq t_0,
\quad x\geq x_0,
\end{align*}
where $B$ is a Brownian sheet.  The coefficients are assumed to be
time-space homogeneous, i.e. they depend on the state variable but not
explicitly on $(t,x)$.  This assumption keeps the notation transparent and is
natural for the Markovian formulation developed below.

A main observation is that the random variable $Y(t,x)$ alone is not a
sufficient state variable.  To restart the system from $(t,x)$, one must also
retain the boundary derivatives along the two coordinate directions.  Thus we
work with the enlarged process
\[
\mathcal{Y}(t,x)=(Y(t,x),D_1Y(t+\cdot,x),D_2Y(t,x+\cdot)),
\]
where $D_1$ and $D_2$ denote differentiation in the first and second
coordinates, respectively.  We prove that this triple satisfies the Markov
property with respect to the natural past sigma-algebra
\[
\mathcal H_{t,x}=\mathcal F_{t,\infty}\vee\mathcal F_{\infty,x}.
\]
In words, conditionally on the horizontal and vertical histories up to
$(t,x)$, the future field has the same law as a newly started It\^o sheet with
initial data given by the current value and the two boundary profiles.

This Markov property differs from the classical one-parameter Markov property
in two essential ways. First, the conditioning is not with respect to the
rectangular past $\mathcal F_{t,x}$ alone, but with respect to the two-sided
past $\mathcal H_{t,x}$, which contains the horizontal and vertical histories
meeting at $(t,x)$. Second, the state variable must contain the boundary
information needed to restart the sheet beyond $(t,x)$. Thus the relevant
Markovian state is not merely the point value $Y(t,x)$, but the current value
together with the two boundary profiles.

The main contribution of this paper is to formulate this enlarged-state Markov
property and to prove it for time-space homogeneous It\^o sheets driven by the
Brownian sheet. In this sense, our result complements the existing Markov-type
theories for Brownian sheets and multiparameter Gaussian fields. Instead of
studying the Brownian sheet itself as a Gaussian random field, we study
stochastic systems driven by the sheet and identify the correct state variable
for their two parameter evolution.

As an additional tool, we recall a Brownian sheet Girsanov theorem and use it to treat drifted sheets.

The paper is organized as follows.  Section~2 introduces time-space
homogeneous It\^o sheets and proves the relevant homogeneity and strong
Markov-type statements.  Section~3 establishes the Markov property for the
enlarged state process.  We then present the Brownian sheet Girsanov theorem
and its density process.

\section{Time-space homogeneous Ito diffusions}

Throughout the paper, $B=\{B(t,x);t,x\geq0\}$ denotes a real-valued Brownian
sheet on a complete probability space $(\Omega,\mathcal F,\P)$.  Thus $B(0,x)=
B(t,0)=0$, $B$ is continuous, centered Gaussian, and
\[
\E[B(t,x)B(s,a)]=(t\wedge s)(x\wedge a),
\qquad s,t,a,x\geq0.
\]
Equivalently, for rectangles
\[
R=(s_1,s_2]\times(a_1,a_2]\subset\R_+^2,
\]
we write
\[
B(R):=B(s_2,a_2)-B(s_1,a_2)-B(s_2,a_1)+B(s_1,a_1),
\]
and $B(R)$ is Gaussian with mean zero and variance equal to the Lebesgue area
$|R|=(s_2-s_1)(a_2-a_1)$.  If two rectangles are disjoint, their Brownian sheet
increments are independent.

For $z=(t,x)$ we put
\[
R_z:=[0,t]\times[0,x],
\qquad
B_z:=B(t,x)=B(R_z),
\qquad
dz:=dt\,dx.
\]
The natural partial order on $\R_+^2$ is denoted by
\[
(s,a)\leq(t,x)
\quad\Longleftrightarrow\quad
s\leq t\text{ and }a\leq x.
\]
The natural Brownian sheet filtration is
\[
\F_{t,x}:=\sigma\{B(s,a);0\leq s\leq t,\ 0\leq a\leq x\}\vee\mathcal N,
\]
where $\mathcal N$ denotes the $\P$-null sets.  We also use the two-sided past
sigma-algebra
\[
\H_{t,x}:=\F_{t,\infty}\vee\F_{\infty,x},
\]
where
\[
\F_{t,\infty}:=\sigma\{B(s,a);0\leq s\leq t,\ a\geq0\}\vee\mathcal N,
\quad
\F_{\infty,x}:=\sigma\{B(s,a);s\geq0,\ 0\leq a\leq x\}\vee\mathcal N.
\]
For a predictable field $\phi$ satisfying the usual square-integrability
condition, the notation
\[
\int_{R_z}\phi(\zeta)B(d\zeta)
\]
means the Walsh integral with respect to the Brownian sheet, viewed as an
orthogonal martingale measure.

In this section we also use the formal differential notation
\[
B(dz)=B(dt\,dx),
\qquad
W(t,x)=\frac{\partial^2B(t,x)}{\partial t\partial x},
\]
where $W$ is space-time white noise.  This notation is only formal; all
stochastic integrals are interpreted in the Walsh sense.

We now introduce the class of time-space homogeneous It\^o sheets considered in the sequel.
Throughout the paper, we assume that
$b,\sigma:\mathbb{R}\to\mathbb{R}$ are globally Lipschitz continuous
and satisfy a linear growth condition, that is, there exist constants
$L,C>0$ such that
\[
|b(y)-b(z)|+|\sigma(y)-\sigma(z)|
\leq L|y-z|,
\qquad y,z\in\mathbb{R},
\]
and
\[
|b(y)|+|\sigma(y)|
\leq C(1+|y|),
\qquad y\in\mathbb{R}.
\]
The boundary data $p$ and $q$ are assumed to belong to
$L^2_{\mathrm{loc}}(\mathbb{R}_+)$.
\begin{definition}
As before let $Y=Y^{t_0,x_0,y,p,q}$ be an Ito sheet of the form
\begin{align} 
&Y(t,x)=y + \int_{t_0}^t p(s)ds +\int_{x_0}^x q(a)da  \nonumber\\
&+ \int_{t_0}^t\int_{x_0}^x b(Y(s,a))dsda+\int_{t_0}^t \int_{x_0}^x \sigma(Y(s,a))B(dsda)\;,\quad t\geq t_0, x \geq x_0.\label{Y}
\end{align}

In differential form the integral equation \eqref{Y} gets the form
\begin{align}\label{Itosheet2}
    &\frac{\partial ^2}{\partial t \partial x} Y(t,x) = b(Y(t,x)) + \sigma(Y(t,x)) \diamond W(t,x);\nonumber\\
    &Y(t_0,x_0)=y,\quad\frac{\partial}{\partial t}Y(t,x_0)=p(t),\quad \frac{\partial}{\partial x}Y(t_0,x)=q(x);\quad t \geq t_0, x\geq x_0.
\end{align}
Next, define the process $\tilde{Y}$ as follows;
\begin{align}
&\tilde{Y}(h,k)=y + \int_0^h D_1 \tilde{Y}(u,0)du + \int_0^k D_2 \tilde{Y}(0,v)dv\\
&+\int_0^h \int_0^k b(\tilde{Y}^{t,x,y,p,q}(u,v))du dv\nonumber\\
&+\int_0^h \int_0^k \sigma(\tilde{Y}^{t,x,y,p,q}(u,v))\diamond \widetilde{W}(u,v)du dv;\quad (s=t+u, a=x+v), \label{Y2}
\end{align}
where $\widetilde{W}(u,v)=W(t+u,x+v); u\geq 0,v\geq 0$. \\

We see by \eqref{Itosheet2} that we can write \eqref{Y} as follows:
\begin{align} \label{Itosheet3}
&Y(t,x)=-y + Y(t,x_0) +Y(t_0,x) \nonumber\\
&+ \int_{t_0}^t\int_{x_0}^x b(Y(s,a))dsda+\int_{t_0}^t \int_{x_0}^x \sigma(Y(s,a))B(dsda)\;,\quad t\geq t_0, x \geq x_0.
\end{align}
\end{definition}

\begin{remark}
Note that we are making an important restriction:
\begin{itemize}
    \item 
In \eqref{Y} we assume that $b(t,x,y)$ and $\sigma(t,x,y)$ do not depend
on $t,x$ but on $y$ only. (This is basically an assumption made for technical simplicity, because the
general case can be reduced to this situation by extending the dimension of the state to include $t$ and $x$.)
\end{itemize}
\end{remark}

The resulting process
$Y(t,x)(\omega)$ will have the following property:

\begin{lemma} \label{hom}
The process $Y(t,x)=Y^{0,0,y,p,q}(t,x)$ is \emph{time-space homogeneous}
in the following sense. For fixed $t,x\geq 0$, define the local boundary
profiles
\begin{align}
    p_{t,x}(u):=D_1Y(t+u,x),\qquad
    q_{t,x}(v):=D_2Y(t,x+v),
    \qquad u,v\geq0.
\end{align}
Then
\begin{align}
Y^{t,x,y,p_Y,q_Y}(t+h,x+k)
\stackrel{\mathrm{law}}{=}
Y^{0,0,y,p_{t,x},q_{t,x}}(h,k),
\qquad h\geq0,\ k\geq0.
\end{align}
In particular, the processes
\[
Y^{t,x,y,p_Y,q_Y}(t+h,x+k)
\quad \mbox{and} \quad
Y^{0,0,y,p_{t,x},q_{t,x}}(h,k),
\qquad h\geq0,\ k\geq0,
\]
have the same distributions, where
\begin{align}
    p_Y(s,x)=D_1Y(s,x)
    \quad \text{and} \quad
    q_Y(t,a)=D_2Y(t,a).
\end{align}
\end{lemma}

\dproof
By \eqref{Y} we have
\begin{align}
&Y^{t,x,y,p_Y,q_Y}(t+h,x+k)\nonumber\\
&=y+\int_{t}^{t+h}D_1Y(s,x)\,ds
+\int_{x}^{x+k}D_2Y(t,a)\,da \nonumber\\
&\quad+\int_{t}^{t+h}\int_{x}^{x+k}
b\bigl(Y^{t,x,y,p_Y,q_Y}(s,a)\bigr)\,ds\,da \nonumber\\
&\quad+\int_{t}^{t+h}\int_{x}^{x+k}
\sigma\bigl(Y^{t,x,y,p_Y,q_Y}(s,a)\bigr)
\diamond W(s,a)\,ds\,da \nonumber\\
&\overset{s=t+u,\ a=x+v}{=}
y+\int_0^h D_1Y(t+u,x)\,du
+\int_0^k D_2Y(t,x+v)\,dv \nonumber\\
&\quad+\int_0^h\int_0^k
b\bigl(Y^{t,x,y,p_Y,q_Y}(t+u,x+v)\bigr)\,du\,dv
\nonumber\\
&\quad+\int_0^h\int_0^k
\sigma\bigl(Y^{t,x,y,p_Y,q_Y}(t+u,x+v)\bigr)
\diamond \widetilde W(u,v)\,du\,dv,
\end{align}
where
\[
\widetilde W(u,v):=W(t+u,x+v),
\qquad u,v\geq0.
\]

Therefore, the process
\begin{align}
    \widetilde Y(h,k)
    :=Y^{t,x,y,p_Y,q_Y}(t+h,x+k)
    \label{tildeY}
\end{align}
satisfies
\begin{align}
\widetilde Y(h,k)
&=y+\int_0^h p_{t,x}(u)\,du
+\int_0^k q_{t,x}(v)\,dv \nonumber\\
&\quad+\int_0^h\int_0^k
b\bigl(\widetilde Y(u,v)\bigr)\,du\,dv \nonumber\\
&\quad+\int_0^h\int_0^k
\sigma\bigl(\widetilde Y(u,v)\bigr)
\diamond\widetilde W(u,v)\,du\,dv.
\end{align}

On the other hand,
\begin{align}
Y^{0,0,y,p_{t,x},q_{t,x}}(h,k)
&=y+\int_0^h p_{t,x}(u)\,du
+\int_0^k q_{t,x}(v)\,dv \nonumber\\
&\quad+\int_0^h\int_0^k
b\bigl(Y^{0,0,y,p_{t,x},q_{t,x}}(u,v)\bigr)\,du\,dv
\nonumber\\
&\quad+\int_0^h\int_0^k
\sigma\bigl(Y^{0,0,y,p_{t,x},q_{t,x}}(u,v)\bigr)
\diamond W(u,v)\,du\,dv.
\end{align}
Since $\widetilde W$ and $W$ have the same $\P$-distribution,
it follows by weak uniqueness of the solution of the stochastic
integral equation \eqref{Y} that
\[
Y^{t,x,y,p_Y,q_Y}(t+h,x+k)
\stackrel{\mathrm{law}}{=}
Y^{0,0,y,p_{t,x},q_{t,x}}(h,k).
\]
Hence Lemma \ref{hom} holds.
\fproof

\section{The Markov property}
From now on we let
$\P^{y,p,q}$ denote the probability law of a given (time-homogeneous) It\^o sheet $\{Y(t,x)\}_{t,x\geq 0}$ when its initial values are $y,p,q$. The expectation with respect to $\P^{y,p,q}$ is denoted by $\E^{y,p,q}[\,\cdot\,]$. Hence we write
\begin{equation}
\E^{y,p,q}[f_1(Y(t_1,x_1))\cdots f_k(Y(t_k, x_k))]=\E[f_1(Y^{y,p,q}(t_1,x_1))\cdots f_k(Y^{y,p,q}(t_k, x_k))]
\end{equation}
for all bounded Borel functions $f_1,\ldots,f_k$ and all times $t_1,x_1,\ldots,t_k\geq 0$; $k=1,2,\ldots$, where $\E=\E_{\P}$ denotes the expectation with respect to the probability law $\P$ of $B(\cdot)$.
Let $\mathcal F_{t,x}$ be the completed $\sigma$-algebra generated by
$\{B(s,a):0\leq s\leq t,\ 0\leq a\leq x\}$.
 Finally we define 
\begin{align}
\mathcal{H}_{t,x}=\mathcal{F}_{t, \infty} \vee \mathcal{F}_{\infty,x} \label{H}
\end{align}
to be the $ \sigma$-algebra generated by $$\{B(s,a); s\leq t, a\geq 0 \}\quad \textnormal{and}  \quad \{B(u,v); u \geq 0, v \leq x\}.$$

For every deterministic $(t,x)$, define
\begin{equation}\label{bhat}
\widehat B_{t,x}(u,v)
:=
B(t+u,x+v)
-B(t+u,x)
-B(t,x+v)
+B(t,x),
\qquad u,v\geq0.
\end{equation}
Indeed, $\widehat B_{t,x}$ is generated by Brownian sheet increments
over rectangles contained in $(t,\infty)\times(x,\infty)$, which are
disjoint from the rectangles generating $\mathcal H_{t,x}$; hence the
independence follows from the independent increment property of the
Brownian sheet.

We now prove that the \emph{triple} 
$$\mathcal{Y}(t,x):=(Y(t,x),D_1 Y(t+\cdot,x), D_2 Y(t,x+\cdot))$$ 
satisfies the \emph{ Markov property}, which basically says that
the future behavior of the process $\mathcal{Y}$ given what has happened up to time-space
$t,x$ is the same as the behavior obtained when starting the process at
$(Y(t,x), D_1 Y(t+\cdot,x), D_2 Y(t,x+\cdot))$.

The notation $D_1Y(t+\cdot,x)$ and $D_2Y(t,x+\cdot)$ is understood in the generalized sense. In particular, these objects need not be represented by ordinary $L^2_{\mathrm{loc}}$ -functions. Whenever they occur through integration, they are identified with the corresponding boundary increments,
\[
\int_0^h D_1Y(t+u,x)du
:=Y(t+h,x)-Y(t,x),
\]
and
\[
\int_0^k D_2Y(t,x+v)dv
:=Y(t,x+k)-Y(t,x).
\]
Equivalently, the state may be formulated directly in terms of these two boundary
increment profiles. The two formulations contain the same boundary information
needed to restart the sheet.

To make the enlarged state precise, let
\[
\mathcal C_0:=\{f\in C(\mathbb R_+):f(0)=0\},
\]
endowed with the topology of uniform convergence on compact intervals
and its Borel $\sigma$-algebra. We identify the generalized boundary
profiles $D_1Y(t+\cdot,x)$ and $D_2Y(t,x+\cdot)$ with their primitives
\[
h\mapsto Y(t+h,x)-Y(t,x),
\qquad
k\mapsto Y(t,x+k)-Y(t,x),
\]
respectively. Thus the enlarged state
\[
\mathbf Y(t,x)
=
\bigl(
Y(t,x),
D_1Y(t+\cdot,x),
D_2Y(t,x+\cdot)
\bigr)
\]
is regarded as an element of the state space
\[
\mathcal E:=\mathbb R\times\mathcal C_0\times\mathcal C_0.
\]
With this identification, $\mathbf Y(t,x)$ is an $\mathcal E$-valued measurable
random variable.

With this notation understood, the precise mathematical formulation is the following:
\begin{theorem}[The Markov property for It\^o sheets]\label{TMarkovProp}
\hfill\break
Let $f$ be a bounded Borel function from $\R$ to $\R$. Then, for
$t,h, x,  k\geq 0$
\begin{equation}\label{mp1}
\E^{y,p,q}[f(Y(t+h,x+k))|\mathcal{H}_{t,x}](\omega)
=\E^{\tilde{Y}(t,x),D_1 \tilde{Y}(t+\cdot,x), D_2 \tilde{Y}(t,x+\cdot))}[f(\tilde{Y}(h,k))].
\end{equation}
where $\tilde{Y}$ is given by \eqref{Y2}.\\
 The right
hand side here means the function $\E^{y,p,q}[f(\tilde{Y}(h,k))]$ evaluated at $y=\tilde{Y}(t,x), p=D_1 \tilde{Y}(t+\cdot,x), q= D_2 \tilde{Y}(t ,x+\cdot).$  Here \(D_1Y(t+\cdot,x)\) and \(D_2Y(t,x+\cdot)\) denote the boundary data of the restarted sheet along the two edges issuing from \((t,x)\).\\
\end{theorem}
\dproof
By the It\^o formula for the Brownian sheet (see Wong and Zakai \cite{WZ} and also Agram et al. \cite{AOPT1})  together with the decomposition of the stochastic integral equation
relative to the point $(t,x)$, we obtain, for $u\ge t$, $v\ge x$,
\begin{align}
    Y(u,v)&=Y(t,x) +\Delta_{u} (Y(t,x))+\Delta_{v}(Y(t,x))\nonumber\\
    &+\int_{t}^{u} \int_{x}^{v} b(Y(s,a))ds da+\int_{t}^{u} \int_{x}^{v} \sigma(Y(s,a)) B(ds da).\label{5.2}
\end{align}
Note that $\Delta_u (Y(t,x))= \int_{t}^{u}\frac{\partial}{\partial s}Y(s,x)ds$ and  $\Delta_v (Y(t,x))= \int_{x}^{v}\frac{\partial}{\partial a}Y(t,a)da.$\\

Therefore we have by uniqueness
$$
Y(u,v)
=
Y^{t,x,Y(t,x),D_1Y(t+\cdot,x),D_2Y(t,x+\cdot)}(u,v)
\quad a.s.
$$
In other words, if we define, for $u\geq t, v\geq x,$
\begin{align}
F(t,x,y,p,q,u,v)
&= Y^{t,x,y,p,q}(u,v) \nonumber\\
&=
y+\int_t^{u} p(s)\,ds+\int_x^{v} q(a)\,da \nonumber\\
&\quad
+\int_t^{u}\int_x^{v}
b\left(Y^{t,x,y,p,q}(s,a)\right)\,ds\,da \nonumber\\
&\quad
+\int_t^{u}\int_x^{v}
\sigma\left(Y^{t,x,y,p,q}(s,a)\right)\,B(ds\,da).
\label{3.5}
\end{align}
we have
\begin{equation}\label{F}
Y(u,v)=F(t,x,Y(t,x),D_1 Y(t+\cdot,x), D_2 Y(t ,x+\cdot),u,v);\;\, u\geq t, v\geq x.
\end{equation}
Fix $(t,x)\in\R_+^2$. Recall that the north-east rectangular
increment $\widehat B_{t,x}$ introduced above \eqref{bhat} is a Brownian sheet
independent of $\mathcal H_{t,x}$.

In the restarted equation below, $p$ and $q$ are understood as the local
boundary profiles, namely
\[
p(r)=D_1Y(t+r,x), \qquad q(a)=D_2Y(t,x+a).
\]
For fixed boundary data \(y,p,q\), we now interpret
\(F(t,x,y,p,q,t+u,x+v)\) as the solution of the restarted equation

\begin{align}
F(t,x,y,p,q,t+u,x+v)
&=
y+\int_0^u p(r)\,dr+\int_0^v q(a)\,da\nonumber \\
&\quad+
\int_0^u\int_0^v
b(F(t,x,y,p,q,t+r,x+a))\,dr\,da\nonumber \\
&\quad+
\int_0^u\int_0^v
\sigma(F(t,x,y,p,q,t+r,x+a))\,\widehat B_{t,x}(dr\,da).
\end{align}

Thus, for fixed \(y,p,q\), the random variable
\[
\omega\mapsto F(t,x,y,p,q,t+u,x+v)(\omega)
\]
is measurable with respect to the sigma-algebra generated by \(\widehat B_{t,x}\).
Since \(\widehat B_{t,x}\) is independent of \(\mathcal H_{t,x}\), we obtain that
\begin{align}\label{ind}
\omega\to F(t,x,y,p,q, t+u,x+v)(\omega) \text{ is independent of } 
\mathcal{H}_{t,x}. 
\end{align}

 Using \eqref{F} we may rewrite the Markov property \eqref{mp1} as follows:
\begin{align} \label{mp2}
&\E[f(F(t,x,Y(t,x),D_1Y(t+\cdot,x),D_2Y(t,x+\cdot),t+h,x+k))|\mathcal H_{t,x}]\nonumber\\
&=\E[f(F(0,0,y,p,q,h,k))]_{y=Y(t,x),\,p=D_1Y(t+\cdot,x),\,q=D_2Y(t,x+\cdot)}.
\end{align}
To prove \eqref{mp2}, put $g(y,p,q,\omega)=f\circ F(t,x,y,p,q,t+h,x+k)(\omega)$.
Then $(y,p,q,\omega)\to g(y,p,q,\omega)$ is measurable. Hence we can approximate $g$ pointwise
boundedly by functions of the form
$$
\sum_{k=1}^m \phi_k(y,p,q)\psi_k(\omega)\;.
$$
Using the properties of conditional expectation and \eqref{ind} we get
\begin{align*}
&\E[g(Y(t,x),D_1Y(t+\cdot,x),D_2Y(t,x+\cdot),\omega)|\mathcal H_{t,x}]\\ &=
   \E\bigg[\lim_{m\to\infty}
\sum_{j=1}^m\phi_j(Y(t,x),D_1Y(t+\cdot,x),D_2Y(t,x+\cdot))\psi_j(\omega)|\mathcal H_{t,x}\bigg] \\
&=\lim_{m\to\infty}
\sum_{j=1}^m\phi_j(Y(t,x),D_1Y(t+\cdot,x),D_2Y(t,x+\cdot))\cdot \E[\psi_j(\omega)|\mathcal H_{t,x}] \\
&=\lim_{m\to\infty}
\sum_{j=1}^m\E[\phi_j(y,p,q)\psi_j(\omega)|\mathcal H_{t,x}]_{y=Y(t,x),\,p=D_1Y(t+\cdot,x),\,q=D_2Y(t,x+\cdot)}\\
&=\E[g(y,p,q,\omega)|\mathcal H_{t,x}]_{y=Y(t,x),\,p=D_1Y(t+\cdot,x),\,q=D_2Y(t,x+\cdot)}\\
&=\E[g(y,p,q,\omega)]_{y=Y(t,x),\,p=D_1Y(t+\cdot,x),\,q=D_2Y(t,x+\cdot)}\;.
\end{align*}
Therefore, since \(\{Y(t,x)\}\) is time-space homogeneous,
\[
\begin{aligned}
&\E\bigl[
f(F(t,x,Y(t,x),D_1Y(t+\cdot,x),D_2Y(t,x+\cdot),t+h,x+k,\omega))
\mid \mathcal H_{t,x}
\bigr] \\
&=
\E\bigl[
f(F(t,x,y,p,q,t+h,x+k,\omega))
\bigr]
\big|_{y=Y(t,x),\,p=D_1Y(t+\cdot,x),\,q=D_2Y(t,x+\cdot)} \\
&=
\E\bigl[
f(F(0,0,y,p,q,h,k))
\bigr]
\big|_{y=Y(t,x),\,p=D_1Y(t+\cdot,x),\,q=D_2Y(t,x+\cdot)}
\end{aligned}
\]
which is \eqref{mp2}.
\fproof
\begin{remark}
\begin{enumerate}
    \item 
Theorem \ref{TMarkovProp} is related in spirit to the sharp Markov property of the Brownian sheet developed by Dalang and Walsh \cite{DalangWal}. Their result concerns the Markovian structure of the Brownian sheet and related random fields relative to suitable conditioning sets. In the present setting, the object of interest is a stochastic system driven by the Brownian sheet. The point of Theorem \ref{TMarkovProp} is to identify the state variable which allows one to restart such a system from a point $(t,x)$. For the It\^o sheets considered here, the point value $Y(t,x)$ alone is not sufficient. The boundary data
$
D_1Y(t+\cdot,x)$ and $D_2Y(t,x+\cdot)
$
along the two edges issuing from $(t,x)$ must also be retained.

    \item The Markov property established in Theorem~\ref{TMarkovProp} is closely relatedto, but different from, the germ Markov property studied by Nualart and Pardoux  \cite{NualartPardoux1994}. Their result concerns the conditional independence of the sigma-fields generated inside and outside an arbitrary domain, given the information on its boundary. In contrast, our result is a dynamic Markov property: the solution on a future rectangle is completely characterized by restarting the stochastic equation from the current point together with the two boundary profiles. 
\end{enumerate}
\end{remark}

\subsection{An infinite-dimensional Markov viewpoint}

We now explain how Theorem~\ref{TMarkovProp} should be interpreted in relation
to the classical Markov theory for multiparameter processes.

It is well known that the Brownian sheet is not a Markov process in the usual
finite-dimensional multiparameter sense. The reason is that the value at a single
point does not contain enough information to determine the future evolution. More
precisely, for a Brownian sheet $B$, the conditional law of
$B(t+h,x+k)$ given the past cannot be determined from $B(t,x)$ alone. Indeed,
one also needs the boundary values along the horizontal and vertical edges
emanating from $(t,x)$.

This is consistent with the viewpoint of Khoshnevisan \cite{Khoshnevisan2002}: the Brownian sheet is not
a multiparameter Markov process in the usual sense, but it admits Markovian
interpretations after enlarging the state space. For example, if one coordinate is
viewed as the time variable, then the Brownian sheet becomes an
infinite-dimensional Feller process with values in a space of continuous functions.

Our result is of the same nature. For the It\^o sheet
\begin{align*}
Y(t,x)
&= y
+ \int_{t_0}^{t} p(s)\,ds
+ \int_{x_0}^{x} q(a)\,da \\
&\quad
+ \int_{t_0}^{t}\int_{x_0}^{x}
b\bigl(Y(s,a)\bigr)\,ds\,da+ \int_{t_0}^{t}\int_{x_0}^{x}
\sigma\bigl(Y(s,a)\bigr)\,B(dsda),
\end{align*}
the scalar value $Y(t,x)$ is not a sufficient Markov state. The correct state is the
enlarged boundary state
\[
\mathbf Y(t,x)
:=
\left(
Y(t,x),
D_1Y(t+\cdot,x),
D_2Y(t,x+\cdot)
\right).
\]
Here $D_1Y(t+\cdot,x)$ and $D_2Y(t,x+\cdot)$ encode the two boundary profiles
needed to restart the sheet beyond the point $(t,x)$.

Let $\mathcal E$ denote a suitable state space of triples
$
(y,p,q),
$
where $y\in\mathbb R$ and $p,q$ are boundary profiles. We endow $\mathcal  E$ with a topology such that convergence
$(y_n,p_n,q_n)\to(y,p,q)$ means $y_n\to y$ and convergence of the
corresponding boundary profiles uniformly on compact intervals.
For a bounded Borel
function $f$, define
\[
P_{h,k}f(y,p,q)
:=
\mathbb E^{y,p,q}
\left[
f(\widetilde Y(h,k))
\right],
\qquad h,k\ge0,
\]
where $\widetilde Y$ denotes the restarted It\^o sheet with initial data
$(y,p,q)$.

Then Theorem~\ref{TMarkovProp} can be written in the Markovian form
\[
\mathbb E
\left[
f(Y(t+h,x+k))
\mid
\mathcal H_{t,x}
\right]
=
P_{h,k}f
\left(
Y(t,x),
D_1Y(t+\cdot,x),
D_2Y(t,x+\cdot)
\right).
\]
Equivalently,
\[
\mathbb E
\left[
f(Y(t+h,x+k))
\mid
\mathcal H_{t,x}
\right]
=
P_{h,k}f(\mathbf Y(t,x)).
\]

Thus the family $\{P_{h,k}\}_{h,k\geq0}$ may be viewed as a transition
family associated with the enlarged infinite-dimensional state space $\mathcal E$.
In this sense, $\mathbf{Y}$ is Markov, although the scalar field $Y(t,x)$ is not.

Under suitable regularity assumptions on $b$ and $\sigma$, for instance Lipschitz
continuity and linear growth, the restarted equation depends continuously on the
initial boundary data $(y,p,q)$. Consequently, if $f$ is bounded and continuous,
then
\[
(y_n,p_n,q_n)\to (y,p,q)
\quad\Longrightarrow\quad
P_{h,k}f(y_n,p_n,q_n)\to P_{h,k}f(y,p,q).
\]
In this case, the operators $P_{h,k}$ are Feller-type operators on the enlarged
state space. Moreover,
\[
\lim_{(h,k)\to(0,0)}
P_{h,k}f(y,p,q)
=
f(y),
\]
whenever the topology on $\mathcal E$ is chosen so that the restarted solution
satisfies
\[
\widetilde Y(h,k)\to y
\qquad\text{as }(h,k)\to(0,0).
\]

Therefore, the Markov property proved here should be understood as an
infinite-dimensional Feller-type Markov property. This is not a finite-dimensional
multiparameter Markov property in the classical sense. Rather, it is a Markov
property obtained after augmenting the state by the two boundary profiles required
to restart the stochastic sheet.

\section{Examples}
The examples in this section illustrate how the Markov property established in Theorem~\ref{TMarkovProp} can be applied to a variety of two-parameter stochastic systems.

\subsection{The Brownian Sheet as a Markov Random Field}\label{ex_BS}
Define
\begin{align}
Y(t,x)=Y^{0,0,0}(t,x)=B(t,x),
\end{align}
starting at $(t_0,x_0)=(0,0)$. After restarting the sheet from $(t,x)$, the corresponding boundary data are denoted by
\[
p(s)=D_1B(t+s,x),\quad s\ge 0,\quad
\textnormal{and} \quad  q(a)=D_2B(t,x+a),\quad a\ge 0.
\]
Equivalently,
\[
\int_0^h p(s)\,ds=B(t+h,x)-B(t,x),
\]
and
\[
\int_0^k q(a)\,da=B(t,x+k)-B(t,x).
\]
Hence, by Theorem \ref{TMarkovProp},
\begin{align*}
&\E[Y(t+h,x+k)| \H_{t,x}]
=
\E^{y,p,q}[\tilde{Y}(h,k)]
\big|_{y=B(t,x),\,p=D_1B(t+\cdot,x),\,q=D_2B(t,x+\cdot)}
\\
&=
\E\Big[\Big(y+\int_0^h p(s)ds +\int_0^k q(a)da +B(h,k)\Big)\Big]
\big|_{y=B(t,x),\,p=D_1B(t+\cdot,x),\,q=D_2B(t,x+\cdot)}
\\
&=
\Big(y+\int_0^h p(s)ds +\int_0^k q(a)da\Big)
\big|_{y=B(t,x),\,p=D_1B(t+\cdot,x),\,q=D_2B(t,x+\cdot)}
+\E[B(h,k)]
\\
&=
B(t,x)+\widetilde B(h,0)-\widetilde B(0,0)
+\widetilde B(0,k)-\widetilde B(0,0)
\\
&=
B(t,x)+B(t+h,x)-B(t,x)+B(t,x+k)-B(t,x)
\\
&=
B(t+h,x)+B(t,x+k)-B(t,x),
\end{align*}
where we have used that $\tilde{B}(u,v)= B(t+u,x+v)$.
In particular, relative to the two-sided sigma-algebra $\mathcal H_{t,x}$, the conditional expectation of $B(t+h,x+k)$ is not determined by the corner value $B(t,x)$ alone, but also involves the horizontal and vertical boundary values $B(t+h,x)$ and $B(t,x+k)$. This illustrates the role of the boundary information in the Markov property considered here.

\subsection{A Nonlinear Multiplicative Noise Sheet}
    The same proof as in the example above shows that if $\sigma:\R \to \R$ is Lipschitz continuous then the process $Y$ defined by
    \begin{align}
    \frac{\partial^2 Y(t,x)}{\partial t \partial x} = \sigma(Y(t,x))  \frac{\partial^2 B(t,x)}{\partial t \partial x}
    \end{align}
    also satisfies the strong martingale property, i.e. for all $t,x,h,k \geq 0$
    \begin{align}
  \E[Y(t+h,x+k)| \H_{t,x}]= Y(t+h,x)+Y(t,x+k) - Y(t,x).
    \end{align}

 \subsection{The two-parameter Langevin equation}
Suppose $Y(t,x)$ is the solution of the Langevin equation
\begin{align}
 dY(z)=\alpha Y(z) dz + \beta B(dz); Y(0)=D_1Y(0)=D_2Y(0)=0,   
\end{align}
where $\alpha$ and $\beta$ are constants.
Then the Markov property gives that
\begin{align*}
    &\E[Y(t+h,x+k)|\H_{t,x}]= \E^{y,p,q}[\widetilde{Y}(h,k)]_{y=Y(t,x), p=D_1Y(t+\cdot,x), q=D_2Y(t,x+\cdot)}\\
    &= \E^{y,p,q}[y+\int_0^h p(s)ds + \int_0^k q(a)da +\int_0^h \int_0^k \alpha \widetilde{Y}(s,a)ds da \\
    &+\beta\int_0^h \int_0^k \widetilde{B}(dsda)]_{y=Y(t,x), p=D_1Y(t+\cdot,x), q=D_2Y(t,x+\cdot)}\\
    &=Y(t,x)+\int_0^h D_1Y(t+s,x)ds + \int_0^k D_2Y(t,x+a)da+\alpha \int_0^h \int_0^k \E^{y,p,q}[\widetilde{Y}(s,a)]dsda\\
    &=Y(t+h,x)+Y(t,x+k)-Y(t,x) + \alpha \int_0^h \int_0^kG(s,a)ds da,
\end{align*}
where $G(t,x):= \E^{y,p,q}[\tilde{Y}(t,x)]= \E^{y,p,q}[Y(t,x)]$ is the solution of the equation
\begin{align*}
    &\frac{\partial^2}{\partial t \partial x} G(t.x)=\alpha G(t,x); \\
    &G(0,0)=Y(0,0), \ D_1 G(\cdot,0)=D_1 Y(\cdot,0), D_2 G(0,\cdot)=D_2 Y(0,\cdot).
\end{align*}

\subsection{Markov Property under an Equivalent Change of Measure}

We now illustrate how the Markov property developed above extends to
Brownian sheets with absolutely continuous drifts through an equivalent
change of measure.

Let
\[
Y(t,x)
=
B(t,x)
-
\int_0^t\int_0^x
\theta(s,a)\,ds\,da,
\qquad
0\le t\le T,\quad 0\le x\le X,
\]
where $B$ is a Brownian sheet and $\theta$ is a jointly measurable random
field.

Following Nualart and Pardoux \cite{NualartPardoux1994}, see also Dalang
and Mueller \cite{DalangMueller2009}, suppose that $\theta$ is adapted to
the filtration
\[
\mathcal G_u
:=
\sigma\{B(t,v):0\le t\le T,\ 0\le v\le u\}
\vee\mathcal N,
\qquad
0\le u\le X,
\]
and satisfies
\[
\E\!\left[
\int_0^T\int_0^X
\theta(t,x)^2\,dt\,dx
\right]
<\infty.
\]

Define
\[
L_u
=
\exp\left(
\int_0^T\int_0^u
\theta(t,x)\,B(dtdx)
-
\frac12
\int_0^T\int_0^u
\theta(t,x)^2\,dt\,dx
\right),
\qquad
0\le u\le X,
\]
and assume that $(L_u)_{0\le u\le X}$ is a martingale. Let
\[
\frac{d\Q}{d\P}
=
L_X .
\]

The Cameron--Martin--Girsanov theorem for the Brownian sheet then implies
that the shifted process
\[
\widetilde B(t,x)
=
B(t,x)
-
\int_0^t\int_0^x
\theta(s,a)\,ds\,da
\]
is a Brownian sheet under $\Q$.

\begin{example}
Since $Y=\widetilde B$ under $\Q$, all the results established previously
for Brownian sheets remain valid under the measure $\Q$.

In particular, by Example~\ref{ex_BS},
\[
\E^{\Q}
\!\left[
Y(t+h,x+k)
\,\middle|\,
\mathcal H_{t,x}
\right]
=
Y(t+h,x)
+
Y(t,x+k)
-
Y(t,x),
\]
for all $0\le t,t+h\le T$ and $0\le x,x+k\le X$.

Therefore, after an equivalent change of measure, a drifted Brownian sheet
becomes a genuine Brownian sheet and consequently enjoys the same strong
martingale and Markov properties.
\end{example}

\section*{Acknowledgments} Nacira Agram gratefully acknowledges financial support from the Swedish Research Council (grant no. 2020-04697). 

\bibliographystyle{plain} 
\bibliography{References}
\end{document}